\documentclass[12pt]{amsart}
\usepackage{amsmath}
\usepackage{amssymb}

\newcommand{\R}{\mathbb{R}}

\newcommand{\Rn}{\mathbb{R}^n}
\newcommand{\Hn}{\mathbb{H}^n}
\renewcommand{\H}{\mathbb{H}}

\renewcommand{\phi}{\varphi}

\newtheorem{theorem}{Theorem}[section]
\newtheorem{lemma}[theorem]{Lemma}

\newtheorem{prop}[theorem]{Proposition}
\theoremstyle{definition}

\newcommand{\Htwo}{\mathbb{H}^2}

\newcommand{\rs}{\widehat{\mathbb{C}}}

\begin{document}
\title{Quasiconformal homogeneity after Gehring and Palka}
\author[P. Bonfert-Taylor, R. Canary and E. Taylor]{Petra Bonfert-Taylor, Richard Canary and Edward Taylor}
\thanks{Richard Canary was partially supported by NSF grant DMS -1306992.}

\maketitle

\section{Introduction}

In a very influential paper \cite{GP76} Gehring and Palka introduced the notions of quasiconformally homogeneous and
uniformly quasiconformally homogeneous  subsets of $\overline{\Rn}$. 
Their motivation was to provide a characterization of quasi-disks,
i.e. domains in $\overline{\Rn}$ which are quasiconformally homeomorphic to the unit disk in $\Rn$. (This paper also
introduced the important concept of the quasihyperbolic metric on a domain in $\overline{\Rn}$.) 
As a generalization, Bonfert-Taylor, Canary,
Martin and Taylor \cite{BCMT05} initiated the study of uniformly quasiconformally homogeneous hyperbolic manifolds.

In this paper, we review the theory of quasiconformally homogeneous subsets of $\overline{\Rn}$ and 
uniformly quasiconformally homogeneous hyperbolic manifolds. We finish
with a discussion of open problems in the theory.

The authors have all benefitted from the inspiration, mathematical and personal, of Fred Gehring's mathematical career.
Petra Bonfert-Taylor was a postdoctoral assistant professor under Fred's supervision, Dick Canary was a colleague
for many years and Ed Taylor was also a postdoctoral assistant professor at the University of Michigan. It is our pleasure
to review the impact of one of Fred Gehring's papers on the field.

\section{Quasiconformally homogeneous subsets of $\overline{\Rn}$}

A domain $D\subset\overline{\Rn}$ is said to be {\em quasiconformally homogeneous} if for all $x,y\in D$, there exists a quasiconformal
map $f:D\to D$ such that $f(x)=y$. In their paper, Gehring and Palka \cite[Lemma 3.2]{GP76} observed that every  domain is quasiconformally
homogenous. We will sketch the proof as variations of the proof will appear later.

\begin{prop}{\rm (Gehring-Palka \cite[Lemma 3.2]{GP76})}
\label{beadsarg}
Every  domain in $\overline{\Rn}$ is  quasiconformally homogeneous.
\end{prop}

\noindent{\em Sketch of proof:}
The proof is based on the following simple lemma which is proved by normalizing and exhibiting
an explicit quasiconformal map.

\begin{lemma}{\rm (Gehring-Palka \cite[Lemma 3.1]{GP76})}
\label{ballest}
If $B$ is a round ball in $\overline{\Rn}$ and  $a,b\in B$,  then there exists
a $K$-quasiconformal map \hbox{$f:\overline{\Rn}\to \overline{\Rn}$} such
that $f(a)=b$, $f|_{\overline{\Rn}-B}=id$, and
$$\log K=d_B(a,b)$$ 
where $d_B$ is the Poincar\'e metric on $B$.
\end{lemma}

If now $D$ is an arbitrary domain in $\mathbb{R}^n$ and $x$ and $y$ lie in $D$, 
then there exists a finite collection $\{B_i\}_{i=0}^n$ of open round balls in $D$ such
that $x\in B_0$, $y\in B_n$ and $B_{i-1}\cap B_{i}$ is non-empty for all $i$. We then choose $x_i\in B_{i-1}\cap B_i$ for all
$i=1,\ldots, n-1$ and set $x_0=x$ and $x_n=y$. Lemma \ref{ballest} implies that for all $i=1,\ldots,n$ there exists a quasiconformal
map  $f_i:D\to D$ such that $f(x_{i-1})=x_i$. So, if $f=f_n\circ\cdots\circ f_1$, then $f:B\to B$ is quasiconformal and
$f(x)=y$.
\qed

\medskip

Thus, it is natural to require that there is a uniform upper bound on the dilatation of the quasiconformal map.
A domain $D\subset\overline{\Rn}$ is said to be {\em $K$-quasiconformally homogeneous} if for all $x,y\in D$, there exists a $K$-quasiconformal
map $f:D\to D$ such that $f(x)=y$.   If $D$ is $K$--quasiconformally homogeneous for
some $K$, we say that it is {\em uniformly quasiconformally homogeneous.}  

A domain $D\subset\overline{\Rn}$ is a {\em $K$-quasi-disk}
if there exists a $K$-quasiconformal homeomorphism $f:\overline{\Rn}\to \overline{\Rn}$ such that $f(\Delta^n)=D$
where $\Delta^n$ is the unit disk in $\Rn$.
Gehring and Palka observed that every $K$-quasi-disk is $K^2$-quasiconformally homogeneous. 
They further show 

\begin{theorem}{\rm (Gehring-Palka \cite[Theorem 5.5]{GP76})}
If $D\subset \overline{\Rn}$ is
a uniformly quasiconformally homogeneous domain with boundary $\partial D$, then either
\begin{enumerate}
\item
$\partial D$ consists of 0,  1 or 2 points,
\item
$\partial D$ is  a non-degenerate continuum, or 
\item
every neighborhood of  every point in $\partial D$ contains infinitely many
components of $\partial D$.
\end{enumerate}
\end{theorem}

Since all simply connected domains in the plane are conformally equivalent to the unit disk, one must further
strengthen these conditions to obtain a characterization of quasi-disks.  One says that a subset $R\subset\overline{\Rn}$
is  {\em $K$-ambiently quasiconformally homogeneous} if there exists $K>0$ 
such that  for all $x,y\in R$, there exists a $K$-quasiconformal
map $f:\overline{\Rn}\to \overline{\Rn}$ such that $f(R)=R$ and $f(x)=y$. 
It is said to be {\em ambiently quasiconformally homogeneous} if it
is $K$-ambiently quasiconformally homogeneous for some $K$. 
Notice that here our subsets need not be domains.

Gehring and Palka constructed
the first non-trivial examples of  ambiently quasiconformally homogeneous domains by observing
that any component of the domain of discontinuity of a convex cocompact subgroup  of ${\rm Isom}_+(\H^{n+1})$ 
is an ambiently quasiconformally homogeneous domain in $\overline{\Rn}=\partial\H^{n+1}$ 
(see \cite[Lemma 4.3]{GP76}).
In particular, by considering Schottky groups, they showed that
there exists an uniformly quasiconformally homogeneous domain in $\overline{\Rn}$ whose complement
is a Cantor set (see \cite[Example 4.4]{GP76}). (In each case, Gehring and Palka only claim uniform quasiconformal homogeneity,
but the proofs they offer immediately establish ambient quasiconformal homogeneity.)

Erkama \cite{ERK77} proved that a Jordan curve in the plane
is ambiently  quasiconformally homogeneous if and only if it is a quasi-circle (i.e. it is the image of the unit circle
under a quasiconformal map of $\rs$).
Brechner and Erkama \cite{BE79} extended this result from Jordan curves to non-degenerate continua. 

Sarvas \cite{SAR85} obtained a characterization of quasi-disks.

\begin{theorem}{\rm (Sarvas \cite{SAR85})}
A domain $D\subset\overline{\R}^2$ is a quasidisk if and only if it is
an ambiently quasiconformally homogeneous Jordan domain.
\end{theorem}

\noindent{\em Sketch of Proof:}
Gehring and Palka \cite{GP76} observed that $K$-quasi-disks are $K^2$-ambiently quasiconformally homogeneous,
since  the unit disk is ambiently conformally homogeneous and the product of
two $K$-quasiconformal maps is $K^2$-quasiconformal.

Now suppose that a Jordan domain $D$ is not a quasi-disk, but is \hbox{$K$-ambiently} quasiconformally homogeneous.
We may assume that $D$ is a bounded domain in $\mathbb C$.
Ahlfors \cite{ahlfors} showed that if $D$ is not a quasi-disk, then  there exists a sequence
$\{(u_n,v_n,w_n)\}$ of triples of distinct points in $\partial D$ so that if $J_n$ and $J_n'$ are the components of 
$\partial D-\{u_n,v_n\}$,
and ${\rm diam}(J_n)\le {\rm diam}(J_n')$, then $w_n\in J_n$ and
$$\lim \frac{|w_n-v_n|}{|u_n-v_n|}=\infty.$$

We may pass to a subsequence, and possibly choose new triples $\{(u_n,v_n,w_n)\}$, so that  
either (a) the open line segment $(u_n,v_n)$ and the round half-disk $D_n$ with
partial boundary $(u_n,v_n)$ (on the same side of $(u_n,v_n)$ as $J_n$)
 are both contained in $D$ for all $n$, or (b) $(u_n,v_n)$ is contained in $\mathbb C-D$ for all $n$.

In case (a), fix a point $a\in D$ and let $f_n:\rs\to\rs$ be a \hbox{$K$-quasiconformal}
map so that $f_n(D)=D$ and $f_n(a)=y_n$ where
$y_n$ is the point in $D_n$ midway ``above'' the midpoint of $(u_,v_n)$. Let $L_n:\rs\to \rs$ be the extension
of an affine map of $\mathbb C$ so that $L_n(u_n)=-1$, $L_n(v_n)=1$ and $L_n(y_n)=\frac{1}{2} i$. Then,
$\{g_n=L_n\circ f_n\}$ is a normal family, since $g_n(D)$ always misses $-1$, $1$ and $\infty$, so, up to subsequence,
it converges to a \hbox{$K$-quasiconformal} map $g:\rs\to\rs$. 
(Notice that $g$ must be non-constant, since 
$g_n(a) = \frac{1}{2} i$ for all $n$ and $\{ {\rm diam}(g_n(D))\}\to\infty$.) 
Now, for each $n$, choose $z_n\in J_n'$ so that $|z_n-v_n|\ge\frac{1}{2} |w_n-v_n|$,
and pass to a subsequence so that the following limits all exist:
$\lim f_n^{-1}(u_n)=\hat u$, $\lim f_n^{-1}(v_n)=\hat v$, $\lim f_n^{-1}(w_n)=\hat w$ and
$\lim f_n^{-1}(z_n)=\hat z$. Then $g(\hat u)=-1$,
$g(\hat  v)=1$ and $g(\hat w)=g(\hat z)=\infty$
(since $\lim \frac{|w_n-v_n|}{|u_n-v_n|}=\infty$). This is impossible since it implies that
$\hat u$, $\hat v$, and $\hat w$ are all distinct, but $\hat z=\hat w$.
However, $\hat w$ and $\hat z$ lie in distinct components of $\partial D-\{\hat u,\hat v\}$. 
The argument to handle case (b) uses similar techniques.
\qed

\medskip

Hjelle \cite{hjelle} showed that the assumption that $D$ is a Jordan domain is necessary in Sarvas' theorem by giving
an example of an ambiently quasiconformally homogeneous simply connected
subset of $\rs$ which is not a quasi-disk. The domain of discontinuity of a purely hyperbolic
degenerate group also provides such an example (see Bonfert-Taylor-Canary-Martin-Taylor-Wolf \cite[Theorem 1.5]{BCMTW10}).

MacManus, N\"{a}kki, and Palka \cite[Theorem 3.1]{MNP98} characterize the possible topological types of ambiently
quasiconformally homogeneous compact subsets of $\rs$.

\begin{theorem}{\rm (MacManus-N\"akki-Palka \cite[Theorem 3.3]{MNP98})} If $R$ is an ambiently quasiconformally homogeneous compact
subset of $\rs$, then either
\begin{enumerate}
\item
$R=\rs$,
\item
$R$ is a finite set of points,
\item
$R$ is a finite union of disjoint quasicircles bounding a domain in $\rs$, or
\item
$R$ is a Cantor set of Hausdorff dimension less than $2$.
\end{enumerate}
\end{theorem}

All sets of type (a), (b) and (c) are ambiently quasiconformally homogeneous, but there is no known characterization of
which Cantor sets are ambiently quasiconformally homogeneous.
However, 
they show that the middle-$\frac{1}{3}$-Cantor set (see \cite[Example 3.5]{MNP98}) and limit sets of  Schottky groups
(see \cite[Theorem 1.2]{MNP99})
are ambiently quasiconformally homogeneous. Therefore,
there exist ambiently quasiconformally homogeneous Cantor sets with any Hausdorff dimension in $(0,2)$
(see also Gong-Martin \cite{gong-martin}).

MacManus, N\"{a}kki, and Palka \cite{MNP99} further define a subset $E\subset \rs$ to 
be {\em uniformly quasiconformally bi-homogeneous}
if there exists $K$ so that if $(a,b),(c,d)\in E\times E^c$ then there exists a $K$-quasiconformal map 
$f:\rs\to\rs$ such that $f(E)=E$,
$f(a)=c$ and $f(b)=d$. They show \cite[Theorem B]{MNP99} that a non-empty compact subset 
of $\rs$ is uniformly quasiconformally 
bi-homogeneous if and only if it either (a) consists of at most two points, (b) is a quasi-circle or 
(c) is an image of the middle-$\frac{1}{3}$-Cantor set under a quasiconformal homeomorphism of $\rs$.
They further show \cite[Theorem D]{MNP99} that a Cantor set $E$ is
uniformly quasiconformally  bi-homogeneous if and only if $E$ is uniformly perfect and $E^c$ is a uniform domain.

Bonfert-Taylor and Taylor \cite[Theorem 1.1]{BT08} show that if $E$ is a Cantor set  in $\rs$ and both 
$E$ and its complement $E^c$ are ambiently
quasiconformally homogeneous, then $E$ is quasiconformally bi-homogenous. Therefore, $E$ must be
uniformly perfect and $E^c$ must be a uniform domain.
Moreover, they exhibit Cantor sets $E$ and $F$ such that 
(a) $E$ is ambiently quasiconformally homogeneous and $E^c$ is not
(see \cite[Example 3.3]{BT08}), and
(b) $F$ is not ambiently quasiconformally homogeneous, but $F^c$ is (see \cite[Example 3.1]{BT08}).

\section{Uniformly quasiconformally homogeneous hyperbolic manifolds}

Inspired by the work of Gehring and Palka \cite{GP76}, Bonfert-Taylor, Canary, Martin and Taylor
\cite{BCMT05} initiated the study of uniformly quasiconformally homogeneous hyperbolic manifolds.
In this paper, {\bf all manifolds will be orientable}.

A (complete) hyperbolic manifold $N=\Hn/\Gamma$ is said to be {\em\hbox{ $K$-quasiconformally} homogeneous} if for
all $x,y\in N$, there exists a \hbox{$K$-quasiconformal} homeomorphism $f:N\to N$ such that $f(x)=y$. It is said to be
{\em uniformly quasiconformally homogeneous} if it is \hbox{$K$-quasiconformally} homogeneous for some $K$. We define
the quasiconformal homogeneity constant of $N$ to be
$$K(N)=\min\{\ K\ |\ N\ {\rm is}\ K{\rm -quasiconformally}\ {\rm homogeneous}\ \}.$$
It is an immediate consequence of compactness theorems for families of $K$-quasiconformal maps that
we may take minimum, rather than simply infimum, in this definition (see \cite[Lemma 2.1]{BCMT05}).

One may use the geometry of quasiconformal homeomorphisms  to obtain some basic restrictions on the geometry
of uniformly quasiconformally homogeneous hyperbolic manifolds. Let $\ell(N)$ denote the infimum of the set of
lengths of homotopically non-trivial closed curves in $N$ and let $d(N)$ denote the supremum of the set of diameters of embedded
hyperbolic balls in $N$.

\begin{theorem}{\rm (Bonfert-Taylor-Canary-Martin-Taylor \cite[Theorem 1.1]{BCMT05})}
\label{basicfacts}
For all $n\ge2$ and $K\ge 1$, there exists $m(n,K)>0$ such that if \hbox{$N=\Hn/\Gamma$}
is a $K$-quasiconformally
homogeneous hyperbolic $n$-manifold, other than $\Hn$, then
\begin{enumerate}
\item
$\ell(N)$ is positive  and $d(N)$ is finite. In particular, 
$$d(N)\le K\ell(N)+2K\log 4.$$
\item
$\ell(N)\ge m(n,K)$, and
\item
every non-trivial element of $\Gamma$ is hyperbolic and the limit set $\Lambda(\Gamma)$ of $\Gamma$ is
all of $\partial\Hn$.
\end{enumerate}
\end{theorem}

\noindent
{\em Sketch of proof:}  Suppose that $x$ lies on a closed homotopically non-trivial curve $\alpha$ of length $\ell$ and that $y$ is
the center of an embedded hyperbolic ball $B$ of radius $r$. Let $f:N\to N$ be a $K$-quasiconformal homeomorphism
such that $f(x)=y$. Since every $K$-quasiconformal homeomorphism is a  $(K,K\log 4)$-quasi-isometry
(see Vuorinen \cite[Theorem 11.2]{vuorinen}) and there exists $z\in\alpha$ such that $f(z)$ is not contained in $B$
(since $f(\alpha)$ is homotopically non-trivial), we see that
$$r\le d(f(x),f(z))\le K d(x,z)+K\log 4\le Kl/2+K\log 4$$
and (1) follows.

Since there is a uniform positive lower bound $d_n$ on $d(N)$ which depends only on $n$,
property (2) follows similarly from the fact, again see \cite[Theorem 11.2]{vuorinen}, that if $f$ is $K$-quasiconformal, then
$$\tanh\left(\frac{d(f(x),f(y))}{2}\right)\le \lambda_n^{1-J}\left(\tanh\left(\frac{d(x,y)}{2}\right)\right)^J$$
where $J=K^{1/(1-n)}$ and $\lambda_n\in [4,2e^{n-1}]$ is the Gr\"otzsch constant. In particular,
one may take $m(n,K)=2\tanh^{-1}\left(\lambda_n^{J-1}\tanh(d_n/2)^{1/J}\right)$.

Since $\ell(N)>0$, $\Gamma$ cannot contain parabolic elements and since $d(N)$ is finite, the limit
set $\Lambda(\Gamma)$ must be all of $\partial\Hn$. Property (3) follows.
\qed

\medskip

Gehring and Palka's proof of Proposition \ref{beadsarg} may be easily adapted to show that every closed hyperbolic manifold is uniformly
quasiconformally homogeneous.  If one keeps careful track of the constants one obtains:

\begin{prop}{\rm (Bonfert-Taylor-Canary-Martin-Taylor \cite[Proposition 2.4]{BCMT05})}
\label{closed}
Every closed hyperbolic $n$-manifold $N$ is uniformly quasiconformally homogeneous. Moreover,
$$K(N)\le \left(e^{\ell(N)\over 4}+1\right)^{2(n-1)\left({4 {\rm diam}(N)\over \ell(N)}+1\right)}$$
where  ${\rm diam}(N)$ is the
diameter of $N$.
\end{prop}

A similar argument  shows that:

\begin{prop}{\rm (Bonfert-Taylor-Canary-Martin-Taylor \cite[Proposition 2.7]{BCMT05})}
\label{coversoforbs}
Every regular cover of a closed hyperbolic orbifold $N$ is uniformly quasiconformally homogeneous. 
\end{prop}

As a consequence of Theorem \ref{basicfacts}  and Proposition \ref{closed} 
one see that a geometrically finite hyperbolic manifold is uniformly
quasiconformally homogeneous if and only if it is closed (see \cite[Corollary 1.2]{BCMT05}).

\medskip

We now discuss rigidity phenomena for uniformly quasiconformally homogeneous hyperbolic
$n$-manifolds where $n\ge 3$. The key tool is McMullen's version of Sullivan's rigidity theorem.

\begin{theorem}{\rm (McMullen \cite[Theorem 2.10]{mcmullen})}
\label{rigidity}
Suppose that \hbox{$N=\Hn/\Gamma$} is a hyperbolic $n$-manifold where $n\ge 3$ and there is an
upper bound on the radius of an embedded hyperbolic ball in $N$. If \hbox{$f:N\to N$} is a quasiconformal
homeomorphism, then $f$ is homotopic to an orientation-preserving isometry.
\end{theorem}

We may then combine Proposition \ref{coversoforbs}, Theorem \ref{basicfacts} and Theorem \ref{rigidity} to show
that:

\begin{theorem}{\rm (Bonfert-Taylor-Canary-Martin-Taylor \cite[Theorem 1.3]{BCMT05})}
\label{QCHcovers}
If $n\ge 3$, a hyperbolic $n$-manifold is uniformly quasiconformally homogeneous if and only if it
is a regular cover of a closed hyperbolic orbifold.
\end{theorem}

\noindent
{\em Sketch of proof:}
We may assume that $N\ne \Hn$, since the result is clearly true when $N=\Hn$.
Proposition \ref{coversoforbs} shows that all regular covers of closed hyperbolic orbifolds are uniformly
quasiconformally homogeneous.

Now suppose that $N=\Hn/\Gamma$ is $K$-quasiconformally homogeneous. 
Recall that if $N=\Hn/\Gamma$ is a hyperbolic manifold and the limit set 
$\Lambda(\Gamma)$ of $\Gamma$ contains more than three points,
then the group $\Theta$ of orientation-preserving isometries of $N$ acts properly discontinuously on $N$.
Therefore,  since $\Lambda(\Gamma)=\partial\Hn$, by Theorem \ref{basicfacts},
$N$ is a regular cover of   the hyperbolic orbifold $N/\Theta$. We will observe that $N/\Theta$ has bounded
diameter, so is closed.

If $x,y\in N$, then there exists a $K$-quasiconformal
map $f:N\to N$ such that $f(x)=y$. McMullen's  Rigidity Theorem \ref{rigidity} implies that there exists an orientation-preserving
isometry $g:N\to N$ which is homotopic to $f$. We may then choose lifts $\tilde f:\Hn\to\Hn$ and \hbox{$\tilde g:\Hn\to\Hn$}, such
that $\phi=\tilde g^{-1}\circ \tilde f$ is $K$-quasiconformal and $\phi$ extends to the identity map on $\Lambda(\Gamma)=\partial \Hn$.

The family of $K$-quasiconformal homeomorphisms of $\Hn$ which restrict to the identity on $\partial \Hn$ is compact, so one
sees immediately that:

\begin{lemma}{\rm (\cite[Lemma 4.1]{BCMT05})}
\label{psibound}
For all $n\ge 2$, there exists an increasing function $\psi_n:(1,\infty)\to (0,\infty)$ such that if
$\phi:\Hn\to \Hn$ is \hbox{$K$-quasiconformal} and extends to the identity on $\partial \Hn$,
then
$$d(x,\phi(x))\le \psi_n(K)$$
for all $x\in \Hn$. Moreover, $\lim_{K\to 1^+}\psi_n(K)=0$.
\end{lemma}

Lemma \ref{psibound} then implies that 
$d(y,g(x))\le \psi_n(K)$. Therefore,
$N/\Theta$ has diameter at most $\psi_n(K)$, so $N$ is a regular cover
of the closed hyperbolic orbifold $N/\Theta$.
\qed\medskip

Since, for all $n$, there is a uniform positive lower bound $r_n$ on the diameter of a closed hyperbolic $n$-orbifold,
and $\lim_{K\to 1^+} \psi_n(K)=0$, we see that there is a uniform lower bound on the quasiconformal homogeneity
constant of  a hyperbolic $n$-manifold other than $\Hn$.

\begin{theorem}{\rm (Bonfert-Taylor-Canary-Martin-Taylor \cite[Theorem 1.4]{BCMT05})}
\label{lowerboundonK}
If $n\ge 3$,  there exists $K_n>1$ such that if $N$ is a uniformly quasiconformally homogeneous 
hyperbolic $n$-manifold other than $\Hn$, then 
$$K(N)\ge K_n.$$
\end{theorem}

One can completely characterize uniformly quasiconformally hyperbolic homogeneous hyperbolic 3-manifolds
with finitely generated fundamental group.

\begin{theorem}{\rm (Bonfert-Taylor-Canary-Martin-Taylor \cite[Theorem 7.1]{BCMT05})}
If $N$ is a non-compact uniformly quasiconformally homogeneous hyperbolic 3-manifold with finitely
generated fundamental group, then there exists a closed hyperbolic 3-manifold $M$ which fibers over the circle
such that $N$ is the cover of $M$ associated to the fiber.
\end{theorem}

\section{Quasiconformally homogeneous surfaces}

 It is natural to ask whether Theorems \ref{QCHcovers} and \ref{lowerboundonK} generalize to
 the setting of hyperbolic surfaces. Since every diffeomorphism of a closed hyperbolic surface is quasiconformal,
 it is clear that McMullen's rigidity Theorem \ref{rigidity} fails for hyperbolic surfaces. Therefore, the proofs outlined
 in the last section do not extend.
 
 We first review attempts to address the following question:
 
 \medskip\noindent
 {\bf Question 1:} {\em Does there exists $K_2>1$ such that if $S$ is a uniformly quasiconformally homogeneous
 surface, then $K(S)\ge K_2$?}
 
 \medskip
 
 Kwakkel and Markovic \cite{KM11} have resolved the question for planar surfaces.
 (Bonfert-Taylor, Canary, Martin, Taylor and Wolf \cite{BCMTW10} had earlier produced a lower
 bound, greater than 1, on the ambient quasiconformal homogeneity constant of
 a planar hyperbolic surface.)
 
 \begin{theorem}{\rm (Kwakkel-Markovic \cite{KM11})}
 There exists $K_{planar}>1$ so that if $S$ is a planar uniformly quasiconformally hyperbolic surface, then
 $$K(S)\ge K_{planar}.$$
 \end{theorem}

Bonfert-Taylor, Bridgeman, Canary and Taylor \cite{BBCT07} exhibited a lower bound on the uniform
quasiconformal homogeneity constant of any hyperelliptic surface.
We recall that a closed hyperbolic surface $S$ of genus $g$ is \emph{hyperelliptic} if it admits a 
conformal involution with  $2g +2$ fixed points.   Hyperelliptic surfaces are known to form a 
$(2g-1)- \text{complex dimensional subvariety}$ of the Moduli space ${M}_g$ of all (isometry classes of) 
closed hyperbolic surfaces of genus $g$.

\begin{theorem}{\rm (Bonfert-Taylor-Bridgeman-Canary-Taylor \cite{BBCT07})}
There exists a constant $K_{\rm{hyp}} > 1$, such that if $S$ is a closed hyperelliptic hyperbolic 
surface, then 
$$K(S) \geq K_{\rm{hyp}}.$$ 
\end{theorem}

\noindent \emph{Sketch of proof:}
If not, there exists a sequence $\{ S_j\}$ of closed hyperelliptic surfaces such that $\lim K(S_j) =1$. We may assume that
$\{S_j\}$ converges geometrically to a surface $S_\infty$. It is easy to check, using compactness results
for families of quasiconformal maps, that $K(S_\infty)=1$, so that $S_\infty=\Htwo.$ In particular, $\lim \ell(S_j)=+\infty$.

We next observe that on any hyperelliptic surface there exist disjoint embedded hyperbolic balls of radius $\frac{l(S_j)}{4}$
about each fixed point of the hyperbolic involution $\phi_j:S_j\to S_j$.
The balls are embedded  by definition of $\ell(S_j)$, so it remains to check that they
are disjoint. If $\beta$ is a geodesic arc joining any two fixed points, then $\beta\cup\phi(\beta)$
is a closed geodesic, so $\beta$ has
length at least $\ell(S_j)/2$. Since $S_j$ has $2g_j+2$ fixed point, where $g_j$ is the genus of $S_j$,
$$(2g_j+2) {\rm area}(B(\ell(S_j)/4))\le {\rm area}(S_j)=4\pi (g_j-1)$$
where $B(\ell(S_j)/4)$ is the ball of radius of $\ell(S_j)/4$ in $\Htwo$.
Therefore,  there is an upper bound on $\ell(S_j)$ which is a contradiction.
\qed

\medskip

More generally, we say that a closed surface $S$ of genus $g$ is \emph{c-fixed point full}, for  $c \in (0, 2]$, if there exists a non-trivial conformal 
automorphism of $S$ having $c(g+1)$ fixed points. In particular, every hyperelliptic surface is $2$-fixed point full.
The argument outlined above easily generalizes to show: 

\begin{theorem}{\rm (Bonfert-Taylor-Bridgeman-Canary-Taylor \cite{BBCT07})}
For each $c \in (0,2]$ there exists $K_c >1$ so that if $S$ is a $c$-fixed point full closed hyperbolic surface,
then $K \geq K_c$.
\end{theorem}

\medskip

One may also modify the question by considering more restrictive forms of quasiconformal homogeneity,
where the arguments of the previous section do apply.
Bonfert-Taylor, Bridgeman, Canary and Taylor \cite{BBCT07} define a 
hyperbolic surface $S$ to be $K$-\emph{strongly quasiconformally homogenous} if for any $x,y \in S$ there exists a $K$-quasiconformal 
homeomorphism $f:S\to S$ such that $f(x)=y$ and $f$ is homotopic to a conformal automorphism of $S$.   
Similarly $S$ is $K$-\emph{extremely quasiconformally homogenous} if for any $x,y \in S$ there exists a $K$-quasiconformal 
homeomorphism $f:S\to S$ such that $f(x)=y$ and $f$ is homotopic to the identity. Gehring and Palka's Lemma \ref{ballest} can be
again used to show that every closed hyperbolic surface is both strongly and extremely quasiconformally homogeneous
  
We denote the strong quasiconformal homogeneity constant of a surface $S$ by $K_{aut}(S)$ and the extreme quasiconformal homogeneity constant of the surface by $K_{0}(S).$    The following relationships are immediate (see \cite[Lemma 6.1]{BBCT07}):

\begin{enumerate}
\item
If $S$ is extremely quasiconformally homogeneous, then $S$ is also strongly quasiconformally homogeneous, and
$$K_0(S) \geq K_{aut}(S) \geq K(S).$$ 
\item
If $S$ is strongly quasiconformally homogeneous, then 
$$K_{aut}(S) \geq K(S).$$
\end{enumerate}

Lemma \ref{psibound} implies that a hyperbolic surface is extremely quasiconformally homogeneous if and only if it is closed
(see \cite[Theorem 6.4]{BBCT07}). The proofs of Theorems \ref{QCHcovers} and \ref{lowerboundonK} immediately generalize
in the strongly quasiconformally homogeneous setting. In order to obtain explicit bounds, we note that Yamada \cite{Y81} proved
that 
$$\tau=\sinh^{-1}\left(\frac{4\cosh^2(\pi/7)-3}{8\cos(\pi/7)+7}\right)\approx 0.13.1467$$
is a lower bound for the diameter of a closed hyperbolic 2-orbifold and that Proposition 6.2 in \cite{BBCT07} provides an explicit
formula for $\psi_2^{-1}$ in terms of the modulus of the Gr\"otsch ring, where $\psi_2$ is the function from Lemma \ref{psibound}.

\begin{theorem}
{\rm (Bonfert-Taylor-Bridgeman-Canary-Taylor \cite[Theorem 6.5]{BBCT07})}
\label{stronglynotsharp}
A hyperbolic surface is strongly quasiconformally homogeneous if and only if it is a regular cover of a closed hyperbolic orbifold.  Moreover,
if $S$ is a strongly quasiconformally homogeneous surface, other than $\Htwo$, then 
$$K_{aut}(S) \geq \psi_2^{-1}(\tau)\approx 1.0595.$$
\end{theorem}

\noindent{\bf Remark:} 
Theorem 6.4 in \cite{BBCT07} shows that if $S$ is a closed hyperbolic surface, then
$$K_0(S)\ge \psi_2^{-1}\left(\sinh^{-1}(\sqrt{2}/3)\right)\approx 1.626.$$

\medskip

Bonfert-Taylor, Martin, Reid and Taylor \cite{BMRT11} obtained a sharp version of Theorem \ref{stronglynotsharp}.
They begin by using the isodiametric inequality to show that the $(2,3,7)$-triangle orbifold, denoted $O_{min}$,
has minimal diameter among all hyperbolic two-orbifolds (see \cite[Proposition 2.2]{BMRT11}).
  
\begin{theorem}
{\rm (Bonfert-Taylor-Martin-Reid-Taylor \cite[Theorem 2.3]{BMRT11})}
\label{Kautsharp}
If $S$ is a strongly quasiconformally homogeneous hyperbolic surface, other than $\Htwo$, then 
$$K_{aut}(S) > \overline{K_{aut}}=\psi_2^{-1}({\rm diam}(O_{min}))\approx 1.36138.$$
Moreover, there exists a sequence $\{ S_j\}$ of regular manifold covers of $O_{min}$ such that
$\lim K(S_j)=\overline{K_{aut}}$.
\end{theorem}

\noindent
{\em Sketch of proof:}
The argument in the previous section establishes that if $S$ is a strongly quasiconformally homogeneous surface,
then 
$$K_{aut}(S)\ge\psi_2^{-1}({\rm diam}\left(O_{min})\right).$$
The proof that the inequality is strict requires a detailed analysis of the extremal maps
for Lemma \ref{psibound}.  The extremal map is unique and 
one demonstrates, via the line element field of this map,
that it can not be realized as a quasiconformal deformation of any non-elementary Fuchsian group.  
The reader is referred to \cite{BMRT11} for details.

One may show that if $\{ S_j\}$ is a sequence of regular manifold covers of $O_{min}$ such that $\lim\ell(S_j)=+\infty$, then
$\lim K(S_j)=\overline{K_{aut}}$. (The existence of such a sequence of covers is guaranteed by the fact that finitely generated Fuchsian groups
are residually finite.) One may assume that $\ell(S_j)>>{\rm diam}(O_{min})$.
Given  $x,y \in S_j$, there exists a conformal automorphism  $g:S_j\to S_j$ such that
$d(g(x),y)\le {\rm diam}(O_{min})$. One may then show that there is a quasiconformal mapping $h:S_j\to S_j$ such that
$h(g(x))=y$, $h$ is the identity off of the ball of radius $\ell(S_j)$ about $g(x)$, and
$$K(h)\le  \psi_2^{-1}(\text{diam} (O_{min}))+ \epsilon(\ell(S_j))$$
where $\epsilon(\ell(S_j))\to 0$ as $\ell(S_j)\to +\infty$. It follows that
$$K(S_j)\le  \psi_2^{-1}(\text{diam} (O_{min}))+ \epsilon(\ell(S_j))$$
so that $\lim K(S_j)=\overline{K_{aut}}$ as desired.
\qed

\medskip

In the setting of closed surfaces it is also natural to restrict the (isotopy class of) the quasiconformal map to
lie in some subgroup of the mapping class group. We recall that the {\em mapping class group}  ${\rm Mod}(S)$ of a 
closed surface $S$ is the set of (isotopy classes) of self-homeomorphisms of $S$. If $H$ is a subgroup of
${\rm Mod}(S)$ we say that $S$ is $H_K$-uniformly quasiconformally homogeneous if for any $x,y\in X$
there exists a \hbox{$K$-quasiconformal} homeomorphism $h:S\to S$ so that $h(x)=y$ and  $[h]\in H$. (See Vlamis \cite{vlamis})
for a more detailed discussion.) The {\em Torelli subgroup} of ${\rm Mod}(S)$ is the subgroup consisting
of homeomorphisms which act trivially on $H_1(S)$. Greenfield \cite{greenfield} and Vlamis \cite{vlamis} have
independently proven:

\begin{theorem}{\rm (Greenfield \cite{greenfield}, Vlamis \cite{vlamis})}
There exists $K_{tor}>1$ such that if $S$ is a closed hyperbolic surface, $H\subset {\rm Mod}(S)$
is the Torelli subgroup and $S$ is $H_K$-quasiconformally homogeneous, then
$$K\ge K_{tor}.$$
\end{theorem}

\noindent
{\bf Remark:} Vlamis \cite{vlamis} obtains similar results for level $r$ congruence subgroups
when $r\ge 3$, finite subgroups  and cyclic subgroups
generated by pure mapping classes.

\bigskip

We now turn our attention to the following question,  which is motivated by Theorem \ref{QCHcovers}.

\medskip\noindent
{\bf Question 2:} {\em Does there exist a geometric characterization of uniformly quasiconformally homogeneous
hyperbolic surfaces?}

\medskip

Proposition \ref{coversoforbs} guarantees that that all covers of closed hyperbolic 2-orbifolds are uniformly quasiconformally
homogeneous. However, one may easily construct a quasiconformal deformation $X$ of a non-compact regular cover $S$ of
a closed hyperbolic 2-orbifold $Q$ which is not itself a regular cover of a closed hyperbolic orbifold
(see Bonfert-Taylor-Canary-Martin-Taylor \cite[Example 5.1]{BCMT05}). Then, since $S$ is uniformly quasiconformally
homogeneous, $X$ is also uniformly quasiconformally homogeneous. (One may guarantee that $X$ is not a regular cover
of a closed hyperbolic 2-orbifold by constructing it to have indiscrete length spectrum, for example.) The key difference
here is that a non-compact regular cover of a closed hyperbolic 2-orbifold has an infinite-dimensional quasiconformal
deformation space, while a non-compact regular cover of a closed hyperbolic $n$-orbifold is quasiconformally rigid if $n\ge 3$.
Therefore, the immediate generalization of Theorem \ref{QCHcovers} does not hold in dimension 2.

One might then optimistically hope that every uniformly quasiconformally homogeneous surface is a quasiconformal
deformation of a regular cover of a closed hyperbolic 2-orbifold. Bonfert-Taylor, Canary, Souto and Taylor \cite{BCST11}
showed that this is not the case.

\begin{theorem}
{\rm (\cite[Theorem 1.1]{BCST11})} 
\label{exotic}
There exists a uniformly quasiconformally homogeneous surface which is not a quasiconformal deformation of the regular cover
of any closed hyperbolic $2$-orbifold.
\end{theorem}

\noindent{\em Sketch of proof:}
Given a connected countable graph $X$, each of whose vertices has valence $d\ge 3$, one may construct a hyperbolic surface $S_X$
by ``thickening up'' $X$.
We first choose a compact hyperbolic surface $F$ with geodesic boundary, such that $F$ is homeomorphic to a sphere with $d$ holes
and each boundary component has length 1. One then obtains $S_X$ by replacing each vertex of $X$ by a copy of $F$ and gluing boundary
components which correspond to the same edge of $X$.

If $\phi$ is an automorphism of $X$, then we may construct a $L$-quasiconformal homeomorphism of $S_X$ which mimics $\phi$, i.e.
it takes a copy of $F$ associated to the vertex $v$ to the copy of $F$ associated to the vertex $\phi(v)$. The quasiconformal dilatation
constant $L$ depends only on our choice of $F$. Since there is a lower bound on the injectivity radius of $S_X$ and the diameter
of each copy of $F$ is constant, we may then use a local version of Lemma \ref{ballest} to show that if $x$ and $y$ lie in a copy of $F$,
then there exists a $M$-quasiconformal homeomorphism of $S_X$ taking $x$ to $y$.  Therefore, if the automorphism group
${\rm Aut}(X)$ of $X$ acts transitively on the vertices of $X$, then $S_X$ will be $LM$-quasiconformally homogeneous.

On the other hand, a regular cover of a closed hyperbolic 2-orbifold is quasi-isometric to the finitely generated Cayley graph of the 
group of deck transformations of the covering map.
Therefore, since quasiconformal maps are quasi-isometries,
any quasiconformal deformation of a regular cover of a closed hyperbolic 2-orbifold, is quasi-isometric to the
Cayley graph of a finitely generated group.

So, in order
to construct a uniformly quasiconformally homogeneous surface which is not  a quasiconformal deformation of a regular cover
of a closed hyperbolic 2-orbifold, it suffices to find a connected, countable graph $X$ whose automorphism group acts
transitively on its set of vertices, which is not quasi-isometric to the Cayley graph of any finitely generated group.
Luckily, Eskin, Fisher and Whyte \cite{FEW12} showed that the
Diestel-Leader graphs (see \cite{DL01}) have automorphism groups that act transitively on their vertices, but are not quasi-isometric to
 the Cayley graph of any finitely generated group. This allows us to complete the proof.
 \qed 
 
\section{Open problems}

Many of the open problems in the field revolve around the motivating questions from the previous section.

\medskip\noindent
{\bf Question 1:} {\em Does there exists $K_2>1$ such that if $S$ is a uniformly quasiconformally homogeneous
surface, then $K(S)\ge K_2$?}

\medskip\noindent
{\bf Question 2:} {\em Does there exist a geometric characterization of uniformly quasiconformally homogeneous
hyperbolic surfaces?}

\medskip

Question 2 is intriguing, but mysterious, so we will focus on Question 1.

It follows from Theorem \ref{basicfacts} and Mumford compactness, that there exists $K_2^g>1$
such that if $S$ is a closed hyperbolic surface of genus $g$, then $K(S)\ge K_2^g$.

\medskip\noindent
{\bf  Problem 1:} {\em Explicitly bound $K_2^g$.}

\medskip

If one had success with the previous problem, one might hope to find a constant which worked
for all closed surfaces.

\medskip\noindent
{\bf Problem 2:} {\em Can one find a bound on $K_2^g$ which is independent of $g$? (i.e. can one find a bound
which works for all closed surfaces?)}

\medskip

It is natural to suspect that $K_2$ (assuming it exists) would be strictly less than $\overline K_{aut}$.

\medskip\noindent
{\bf Problem 3:} {\em
Construct a uniformly quasiconformally homogeneous surface $S$ such
that $K(S)<\overline K_{aut}$.}

\medskip

Returning to higher dimensions, one would like to explicitly bound $K_n$.

\medskip\noindent
{\bf Problem 4:} {\em Explicitly bound $K_n$ for all $n\ge 3$.}

\medskip

Manojlovi\'c-Vuorinen \cite{MV} and Vuorinen-Zhang \cite{VZ} have obtained analogues of Lemma \ref{psibound} for
all $n$, with explicit bounds on the resulting functions $\psi_n$. One may use estimates of Adeboye-Wei \cite{adeboye-wei}
to obtain an explicit lower bound, for all $n$, on the diameter of any hyperbolic \hbox{$n$-orbifold}. Since Theorem \ref{basicfacts}
and McMullen's Rigidity Theorem (Theorem \ref{rigidity}) imply that every quasiconformal automorphism of
a uniformly quasiconformally homogeneous hyperbolic manifold is homotopic to an isometry, one should be able to
follow the proof of Theorem \ref{stronglynotsharp} to produce a lower bound on $K_n$. 

\medskip

One might further hope to use the technique of proof of Theorem \ref{Kautsharp} to produce a
sharp lower bound on $K_n$ in higher dimensions. 

In dimension 3, there is a natural candidate for the minimal diameter hyperbolic orbifold.
Gehring and Martin \cite{GM94}  demonstrated that, amongst all Kleinian groups containing a torsion element of order $p \geq 4$,
the unique Kleinian group of minimal co-volume is a $\mathbb Z_2$-extension of the orientation preserving index $2$ subgroup 
generated by reflections in the sides of the hyperbolic tetrahedron with Coxeter diagram $3-5-3$.
Gaven Martin has conjectured:

\medskip\noindent
{\bf Conjecture:} (Martin) {\em The minimum diameter hyperbolic orbifold is the  $\mathbb Z_2$-extension of the orientation-preserving index $2$ subgroup of the $3-5-3$ Coxeter group described above}.

\medskip

This leads us to:

\medskip\noindent
{\bf Problem 5:} {\em Determine $K_3$ precisely. Is it true that if $N$ is a quasiconformally homogeneous hyperbolic 3-manifold,
other than $\mathbb H^3$, then $K(N)>K_3$?} 

\medskip
In order to adapt the proof of Theorem \ref{Kautsharp}, one would
also have to investigate an analogue of the Teichm\"uller extremal map in the three-dimensional setting,
that is, a quasiconformal mapping of minimal distortion that maps $\Delta^n$ to itself, moves the origin to a point $x\neq 0$ 
and extends to the identity on $\partial\Delta^n$.  
However, in dimensions three and above, distortion can be measured in terms of several dilatation functions, 
e.g. the trace dilatation,
the outer and inner dilatations, linear dilatatons, and mean dilatations (in the setting of mappings with finite distortion). 
The existence and uniqueness properties for such problems  depend on the choice of dilatation (see Fehlmann \cite{Feh87}). 
For instance, it is known (see K\"uhnau \cite{Kuh62}) that the extremal problem with boundary data for 
the box problem of Gr\"{o}tsch admits
no unique solution if the dilatation is measured in terms of the inner and outer dilatation functions.  
Progress to date on this problem includes  foundational work by Gehring and Vaisala \cite{GV62} on the extremal problem in the absence
of boundary conditions and work by Astala, Iwaneic, Martin and Onninen \cite{AIMO05}  on the extremal problem for functions of finite distortion.

\end{document}